\documentclass[12pt,smallextended]{svjour3}
\smartqed
\usepackage{lineno}

\journalname{Journal of XXX}
\usepackage{amssymb,amsmath,amsfonts}
\usepackage{hyperref}
\hypersetup{
    colorlinks=true,
    linkcolor=blue,
    filecolor=blue,
    urlcolor=blue,
    citecolor=cyan,
    pdfstartview=FitW,
    bookmarks=true,
}
\usepackage{array}
\usepackage{epsfig}
\usepackage{mathptmx}
\usepackage{diagbox}
\usepackage{graphicx}
\usepackage{color}
\usepackage{caption}
\usepackage{wrapfig}
\usepackage{graphics}
\usepackage{enumerate}
\usepackage{amsxtra}
\usepackage{latexsym}
\usepackage{multirow}
\usepackage{slashbox}
\usepackage{subfigure}
\usepackage{epstopdf}
\usepackage{pdfpages}
\usepackage{algorithm}
\usepackage{algorithmic}
\usepackage{url}
\newcommand{\thmlist}{
\begin{list}{Step 1}
{\setlength{\leftmargin}{0.6 in}\setlength{\labelwidth} {0.5 in}}}

\newcommand{\alglist}{
\begin{list}{Step 1}
{\setlength{\leftmargin}{1.1 in} \setlength{\labelwidth}{1.0 in}}}

 \renewcommand{\proof} {\noindent {\bf Proof.} \quad}

\renewcommand{\subtitle}[1]{\color{blue}}

\def\red#1{\color{red}{#1}\color{black}}

\begin{document}



\title{Continuation Newton methods with the residual trust-region
time-stepping scheme for nonlinear equations}
\titlerunning{Continuation Newton methods for nonlinear equations}
\author{Xin-long Luo \textsuperscript{$\ast$} \and Hang Xiao \and Jia-hui Lv}
\authorrunning{Luo, Xiao and Lv}

\institute{
     Xin-long Luo, Corresponding author
     \at
     School of Artificial Intelligence,
     Beijing University of Posts and Telecommunications, P. O. Box 101,
     Xitucheng Road  No. 10, Haidian District, 100876, Beijing China \\
     \email{luoxinlong@bupt.edu.cn}
     \and
     Hang Xiao
     \at
     School of Artificial Intelligence,
     Beijing University of Posts and Telecommunications, P. O. Box 101,
     Xitucheng Road  No. 10, Haidian District, 100876, Beijing China \\
     \email{xiaohang0210@bupt.edu.cn}
     \and
     Jia-hui Lv
     \at
     School of Artificial Intelligence,
     Beijing University of Posts and Telecommunications, P. O. Box 101,
     Xitucheng Road  No. 10, Haidian District, 100876, Beijing China \\
     \email{jhlv@bupt.edu.cn}
}

\date{Received: date / Accepted: date}

\maketitle

\begin{abstract}
For nonlinear equations, the homotopy methods (continuation methods) are
popular in engineering fields since their convergence regions are large
and they are quite reliable to find a solution. The disadvantage of the
classical homotopy methods is that their computational time is heavy since
they need to solve many auxiliary nonlinear systems during the intermediate
continuation processes. In order to overcome this shortcoming, we consider the
special explicit continuation Newton method with the residual trust-region
time-stepping scheme for this problem. According to our numerical experiments,
the new method is more robust and faster to find the required solution of the
real-world problem than the traditional optimization method  (the built-in subroutine
fsolve.m of the MATLAB environment) and the homotopy continuation methods
(HOMPACK90 and NAClab). Furthermore, we analyze the global convergence and the local
superlinear convergence of the new method.
\end{abstract}




\section{Introduction}

In engineering fields, we often need to solve the equilibrium state of the
differential equation \cite{Liao2012,MM1990,Robertson1966,VL1994}
as follows:
\begin{align}
  \frac{dx}{dt} = F(x), \; x(t_0) = x_{0}.  \label{ODE}
\end{align}
That is to say, it requires to solve the following system of nonlinear equations:
\begin{align}
  F(x) = 0,   \label{NLEFX}
\end{align}
where $F: \; \Re^{n} \to \Re^{n}$ is a vector function. For the nonlinear
system \eqref{NLEFX}, there are many popular traditional optimization methods
\cite{CGT2000,DS2009,Higham1999,NW1999} and the classical homotopy continuation
methods \cite{AG2003,Doedel2007,OR2000,WSMMW1997} to solve it.

\vskip 2mm

For the traditional optimization methods such as the trust-region methods 
\cite{More1978,Yuan1998,Yuan2015} and the line search methods 
\cite{Kelley2003,Kelley2018}, the solution $x^{\ast}$ of the nonlinear system
\eqref{NLEFX} is found via solving the following equivalent nonlinear least-squares
problem
\begin{align}
    \min_{x \in \Re^{n}} f(x) = \frac{1}{2} \|F(x)\|^2,         \label{UOPTF}
\end{align}
where $\|\cdot\|$ denotes the Euclidean vector norm or its induced
matrix norm. Generally speaking, the traditional optimization methods based on
the merit function \eqref{UOPTF} are efficient for the large-scale problems since
they have the local superlinear convergence near the solution $x^{\ast}$
\cite{CGT2000,NW1999}.

\vskip 2mm

However, the line search methods and the trust-region methods are apt to
stagnate at a local minimum point $x^{\ast}$ of problem \eqref{UOPTF}, when
the Jacobian matrix $J(x^{\ast})$ of $F(x^{\ast})$ is singular or nearly
singular, where $J(x) = {\partial F(x)}/{\partial x}$ (see p. 304, \cite{NW1999}).
Furthermore, the termination condition
\begin{align}
   \|\nabla f(x_k)\| = \|J(x_k)^{T} F(x_k)\| < \epsilon, \label{TUNOPT}
\end{align}
may lead these methods to early stop far away from the local minimum $x^{\ast}$.
It can be illustrated as follows. We consider
\begin{align}
  F(x) = Ax = 0,  \hskip 2mm  A =
     \begin{bmatrix}
     1 & 0 \\
     0 & 10^{-6}
    \end{bmatrix}. \label{TINYEXA}
\end{align}
It is not difficult to know that the linear system \eqref{TINYEXA} has a unique
solution $x^{\ast} = (0, \, 0)$. If we set $\epsilon = 10^{-6}$, the traditional
optimization methods will early stop far away from
$x^{\ast}$ provided that $x_{k} = (0, \; c), \; c < 10^{6}$.

\vskip 2mm

For the classical homotopy methods, the solution $x^{\ast}$ of the nonlinear
system \eqref{NLEFX} is found via constructing the following homotopy function
\begin{align}
  H(x, \, \lambda) = (1- \lambda)G(x) + \lambda F(x), \label{HOMOFUN}
\end{align}
and attempting to trace an implicitly defined curve $\lambda(t) \in H^{-1}(0)$
from the starting point $(x_{0}, \, 0)$ to a solution $(x^{\ast}, \, 1)$ by the
predictor-corrector methods \cite{AG2003,Doedel2007}, where the zero point of
the artificial smooth function $G(x)$ is known. Generally speaking, the homotopy
continuation methods are more reliable than the merit-function methods and they
are very popular in engineering fields \cite{Liao2012}. The disadvantage of the
classical homotopy methods is that they require significantly more function and
derivative evaluations, and linear algebra operations than the merit-function
methods since they need to solve many auxiliary nonlinear systems during the
intermediate continuation processes.

\vskip 2mm

In order to overcome this shortcoming of the traditional homotopy methods, we
consider the special continuation method based on the following Newton flow
\cite{AS2015,Branin1972,Davidenko1953,Tanabe1979}
\begin{align}
  \frac{dx(t)}{dt} = - J(x)^{-1}F(x), \;  x(t_0) = x_0,  \label{NEWTONFLOW}
\end{align}
and construct a special ODE method with the new time-stepping scheme based on
the trust-region updating strategy to follow the trajectory of the Newton flow
\eqref{NEWTONFLOW}. Consequently, we obtain its steady-state solution $x^{\ast}$, 
i.e. the required solution $x^{\ast}$ of the nonlinear system \eqref{NLEFX}.

\vskip 2mm

The rest of this article is organized as follows. In the next section, we
consider the explicit continuation Newton method with the trust-region updating 
strategy for nonlinear equations. In section 3, we prove the global convergence 
and the local superlinear convergence of the new method under some standard 
assumptions. In section 4, some promising numerical results of the new method 
are also reported, in comparison to the traditional trust-region method 
(the built-in subroutine fsolve.m of the MATLAB environment
\cite{MATLAB,More1978}) and the classical homotopy continuation methods
(HOMOPACK90 \cite{WSMMW1997} and NAClab \cite{LLT2008,ZL2013,Zeng2019}).
Finally, some conclusions and the future work are discussed in section 5.
Throughout this article, we assume that $F(\cdot)$ exists the
zero point $x^{\ast}$.

\section{Continuation Newton methods}

In this section, based on the trust-region updating strategy, we construct a new
time-stepping scheme for the continuation Newton method to follow the
trajectory of the Newton flow and obtain its steady-state solution $x^{\ast}$.

\subsection{The continuous Newton flow} \label{SUBSECNF}

If we consider the damped Newton method with the line search strategy for
the nonlinear system \eqref{NLEFX} \cite{Kelley2003,NW1999}, we have
\begin{align}
  x_{k+1} = x_{k} - \alpha_{k} J(x_{k})^{-1} F(x_{k}). \label{NEWTON}
\end{align}
We denote $o(\alpha)$ as the higher-order infinitesimal of $\alpha$, that is to
say,
\begin{align}
   \lim_{\alpha \to 0} \frac{o(\alpha)}{\alpha} = 0. \nonumber
\end{align}
In equation \eqref{NEWTON}, if we let  $x_{k} = x(t_{k})$,  
$x_{k+1} = x(t_{k} + \alpha_{k}) + o(\alpha_{k})$ and $\alpha_{k} \to 0$, 
we obtain the continuous Newton flow \eqref{NEWTONFLOW}. Actually, if we apply
an iteration with the explicit Euler method \cite{HNW1993,SGT2003} to the
continuous Newton flow \eqref{NEWTONFLOW}, we also obtain the damped Newton
method \eqref{NEWTON}. Since the Jacobian matrix $J(x)$ may be singular,  we
reformulate the continuous Newton flow \eqref{NEWTONFLOW} as the more general
formula:
\begin{align}
  -J(x)\frac{dx(t)}{dt} = F(x), \hskip 2mm  x(t_0) = x_0. \label{DAEFLOW}
\end{align}

\vskip 2mm

The continuous Newton flow \eqref{DAEFLOW} is an old method and can be backtracked
to Davidenko's work \cite{Davidenko1953} in 1953. After that, it was investigated
by Branin \cite{Branin1972}, Deuflhard et al \cite{PHP1975}, Tanabe \cite{Tanabe1979}
and Kalaba et al \cite {KZHH1977} in 1970s, and applied to nonlinear boundary
problems by Axelsson and Sysala \cite{AS2015} recently. The continuous and even 
growing interest in this method originates from its some nice properties. One of 
them is that the solution $x(t)$ of the continuous Newton flow converges to 
the steady-state solution $x^{\ast}$ from any initial point $x_{0}$, as described 
by the following property \ref{PRODAEFLOW}.

\vskip 2mm 

\begin{property} (Branin \cite{Branin1972} and Tanabe \cite{Tanabe1979})
\label{PRODAEFLOW} Assume that $x(t)$ is the solution of the continuous Newton
flow \eqref{DAEFLOW}, then $f(x(t)) = \|F(x)\|^{2}$ converges to zero when
$t \to \infty$. That is to say, for every limit point $x^{\ast}$ of $x(t)$, it is
also a solution of the nonlinear system \eqref{NLEFX}. Furthermore, every
element $F^{i}(x)$ of $F(x)$ has the same convergence rate $e^{-t}$
and $x(t)$ can not converge to the solution $x^{\ast}$ of the nonlinear system
\eqref{NLEFX} on the finite interval when the initial point $x_{0}$ is not a
solution of the nonlinear system \eqref{NLEFX}.
\end{property}
\proof Assume that $x(t)$ is the solution of the continuous Newton flow
\eqref{DAEFLOW}, then we have
\begin{align}
    \frac{d}{dt} \left(e^{t}F(x)\right) = e^{t} J(x) \frac{dx(t)}{dt}
    + e^{t} F(x) = 0. \nonumber
\end{align}
Consequently, we obtain
\begin{align}
     F(x(t)) = F(x_0)e^{-t}. \label{FUNPAR}
\end{align}
From equation \eqref{FUNPAR}, it is not difficult to know that every element
$F^{i}(x)$ of $F(x)$ converges to zero with the linear convergence
rate $e^{-t}$ when $t \to \infty$. Thus, if the solution $x(t)$ of the continuous
Newton flow \eqref{DAEFLOW} belongs to a compact set, it has a limit point
$x^{\ast}$ when $t \to \infty$, and this limit point $x^{\ast}$ is
also a solution of the nonlinear system \eqref{NLEFX}.

\vskip 2mm

If we assume that the solution $x(t)$ of the continuous Newton flow
\eqref{DAEFLOW} converges to the solution $x^{\ast}$ of the nonlinear system
\eqref{NLEFX} on the finite interval $(0, \, T]$, from equation \eqref{FUNPAR},
we have
\begin{align}
     F(x^{\ast}) = F(x_{0}) e^{-T}. \label{FLIMT}
\end{align}
Since $x^{\ast}$ is a solution of the nonlinear system \eqref{NLEFX}, we have
$F(x^{\ast}) = 0$. By substituting it into equation \eqref{FLIMT}, we
obtain
\begin{align}
     F(x_{0}) = 0. \nonumber
\end{align}
Thus, it contradicts the assumption that $x_{0}$ is not a solution of the nonlinear
system \eqref{NLEFX}. Consequently, the solution $x(t)$ of the continuous Newton flow
\eqref{DAEFLOW} can not converge to the solution $x^{\ast}$ of the nonlinear system
\eqref{NLEFX} on the finite interval. \qed

\vskip 2mm

\begin{remark}
The inverse $J(x)^{-1}$ of the Jacobian matrix $J(x)$ can be regarded as the
preconditioner of $F(x)$ such that the solution elements $x^{i}(t)\, (i = 1, \, 2, \ldots, n)$
of the continuous Newton flow \eqref{NEWTONFLOW} have the roughly same convergence
rates and it mitigates the stiff property of the ODE \eqref{NEWTONFLOW}
(the definition of the stiff problem can be found in \cite{HW1996} and references
therein). This property is very useful since it makes us adopt the explicit ODE
method to follow the trajectory of the Newton flow.
\end{remark}

\vskip 2mm

Actually, if we consider $F(x) = Ax$, from the ODE \eqref{DAEFLOW},
we have
\begin{align}
  A \frac{dx}{dt} = - A x, \; x(0) = x_{0}.  \label{LINAPP}
\end{align}
By integrating the linear ODE \eqref{LINAPP}, we obtain
\begin{align}
   x(t) = e^{-t} x_{0}.       \label{LINSOL}
\end{align}
From equation \eqref{LINSOL}, we know that the solution $x(t)$ of the ODE
\eqref{LINAPP} linearly converges to zero with the same rate $e^{-t}$
when $t$ tends to infinity.

\subsection{Continuation Newton methods}

From subsection \ref{SUBSECNF}, we know that the solution $x(t)$ of the
continuous Newton flow \eqref{DAEFLOW} has the nice global convergence property.
On the other hand, when the Jacobian matrix $J(x)$ is singular or nearly singular,
the ODE \eqref{DAEFLOW} is the system of differential-algebraic equations (DAEs)
and its trajectory can not be efficiently followed by the general ODE method such as
the backward differentiation formulas (the built-in subroutine ode15s.m of the
MATLAB environment \cite{AP1998,BCP1996,HW1996,MATLAB,SGT2003}). Thus, we need to
construct the special method to handle this problem. Furthermore, we expect that 
the new method has the global convergence as the homotopy continuation methods and
the fast convergence rate near the solution $x^{\ast}$ as the merit-function
methods. In order to achieve these two aims, we construct the special continuous
Newton method with the new step size $\alpha_{k} = \Delta t_{k}/(1+\Delta t_{k})$
and the time step $\Delta t_{k}$ is adaptively adjusted by the trust-region
updating strategy for problem \eqref{DAEFLOW}.

\vskip 2mm

Firstly, we apply the implicit Euler method to the continuous Newton flow
\eqref{DAEFLOW} \cite{AP1998,BCP1996}, then we obtain
\begin{align}
  J(x_{k+1}) {(x_{k+1}-x_{k})}/{\Delta t_k} =  -F(x_{k+1}).  \label{IMEDAE}
\end{align}
The scheme \eqref{IMEDAE} is an implicit method. Thus, it needs to solve a system
of nonlinear equations at every iteration. To avoid solving the system of nonlinear
equations, we replace $J(x_{k+1})$ with $J(x_{k})$ and substitute $F(x_{k+1})$
with its linear approximation $F(x_k)+J(x_k)(x_{k+1}-x_{k})$ in
equation \eqref{IMEDAE}. Thus, we obtain the continuation Newton method
as follows:
\begin{align}
  J(x_k) s_{k} & = -  ({\Delta t_k}/{(1+\Delta t_k)})F(x_k), \label{SMEDAE} \\
    x_{k+1} & = x_{k} + s_k.                        \label{XK1}
\end{align}

\vskip 2mm

\begin{remark}
The explicit continuation Newton method \eqref{SMEDAE}-\eqref{XK1} is similar 
to the damped Newton method \eqref{NEWTON} if we let $\alpha_{k} =
\Delta t_k/(1+\Delta t_k)$ in equation \eqref{SMEDAE}. However, from the view 
of the ODE method, they are different. The damped Newton method \eqref{NEWTON} 
is obtained by the explicit Euler method applied to the continuous Newton flow 
\eqref{DAEFLOW}, and its time step $\alpha_k$ is restricted by the numerical 
stability \cite{HW1996,SGT2003}. That is to say, the large time step 
$\alpha_{k}$ can not be adopted in the steady-state phase. The explicit 
continuation Newton method \eqref{SMEDAE}-\eqref{XK1} is obtained by the implicit 
Euler method and its linear approximation applied to the continuous Newton flow 
\eqref{DAEFLOW}, and its time step $\Delta t_k$ is not restricted by the numerical 
stability for the linear test equation $dx/dt = - \lambda x, \; \lambda > 0$. 
Therefore, the large time step $\Delta t_{k}$ can be adopted in the steady-state 
phase for the explicit continuation Newton method \eqref{SMEDAE}-\eqref{XK1}, 
and it mimics the Newton method near the steady-state solution $x^{\ast}$ such that 
it has the fast local convergence rate. The most of all, the new time step 
$\alpha_{k} = \Delta t_{k}/(\Delta t_{k} + 1)$ is favourable to adopt the trust-region 
updating strategy for adaptively adjusting the time step $\Delta t_{k}$ such that 
the explicit continuation Newton method \eqref{SMEDAE}-\eqref{XK1} accurately follows
the trajectory of the continuous Newton flow in the transient-state phase and achieves 
the fast convergence rate near the steady-sate solution $x^{\ast}$.
\end{remark}

For the real-world problem, the Jacobian matrix $J(x)$ may be singular, which
arises from the physical property. For example, for the chemical kinetic reaction
problem \eqref{ODE}, the elements of $x(t)$ represent the reaction concentrations
and they must satisfy the linear conservation law \cite{Logan1996}. A system is
called to satisfy the linear conservation law (\cite{Shampine1998},
or p. 35, \cite{SGT2003}), if there is a constant vector $c \neq 0$ such that
\begin{align}
   c^{T}x(t) = c^{T}x(0) \label{LCL}
\end{align}
holds for all $t \ge 0$. If there exists a constant vector $c$ such
that
\begin{align}
   c^{T}F(x) = 0, \quad \forall x \in \Re^{n},    \label{CFX}
\end{align}
we have
\begin{align}
  c^{T} J(x) = 0,  \quad \forall x \in \Re^{n}. \label{SINJAC}
\end{align}
From equation \eqref{SINJAC}, we know that the Jacobian matrix $J(x)$ is singular.
For this case, the solution $x(t)$ of the ODE \eqref{ODE} satisfies the linear
conservation law \eqref{LCL}.

\vskip 2mm

For the isolated singularity of the Jacobian matrix $J(x)$, there are some
efficient approaches to handle this problem \cite{Griewank1985}. Here, since
the singularity set of the Jacobian matrix $J(x)$ may be connected, we adopt the
regularization technique \cite{Hansen1994,KK1998} to modify the explicit
continuation Newton method \eqref{SMEDAE}-\eqref{XK1} as follows:
\begin{align}
  \left(\mu_{k} I - J(x_k) \right) s_{k}^{P} & =   F(x_k), \label{RSMEDAE} \\
    s_{k} & = ({\Delta t_k}/{(1+\Delta t_k)}) s_{k}^{P},  \label{RSCORDAE} \\
    x_{k+1} & = x_{k} + s_k,                  \label{XKSK}
\end{align}
where $\mu_{k}$ is a small positive number. In order to achieve the fast
convergence rate near the solution $x^{\ast}$, the regularization continuation
Newton method \eqref{RSMEDAE}-\eqref{XKSK} is required to approximate the
Newton method $x_{k+1} = x_{k}- J(x_{k})^{-1}F(x_{k})$  near the solution $x^{\ast}$
\cite{DM1974}. Thus, we select the regularization parameter $\mu_{k}$ as follows:
\begin{align}
   \mu_{k} =
     \begin{cases}
        c_{\epsilon}, \; \text{if} \; \Delta t_{k} \le 1/c_{\epsilon},  \\
        1/\Delta t_{k}, \; \text{others},
     \end{cases}    \label{MUSW}
\end{align}
where $c_{\epsilon}$ is a small positive constant such as $c_{\epsilon} = 10^{-6}$
in practice.

\vskip 2mm

It is not difficult to verify that the regularization continuation Newton method
\eqref{RSMEDAE}-\eqref{XKSK} preserves the linear conservation law \eqref{LCL}
if it exists a constant vector $c \in \Re^{n}$ such that $c^{T}F(x) = 0,
\; \forall x\in \Re^{n}$. Actually, from $c^{T}F(x) = 0$, we have $c^{T}J(x) = 0$.
Therefore, from equations \eqref{RSMEDAE}-\eqref{XKSK}, we obtain
\begin{align}
  c^{T}x_{k+1}  = c^{T}x_{k} + c^{T}s_k
   =  c^{T}x_k + \frac{1}{\mu_{k}}c^{T}\left(\frac{\Delta t_k}{1+\Delta t_k}F(x_k)
  +J(x_k)s_k \right)  = c^{T}x_{k}. \label{CTX}
\end{align}
That is to say, the regularization continuation Newton method \eqref{RSMEDAE}-\eqref{XKSK}
preserves the linear conservation law \eqref{LCL}.

\subsection{The residual trust-region time-stepping scheme}

Another issue is how to adaptively adjust the time-stepping size $\Delta t_k$
at every iteration. A popular way to control the time-stepping size is based on
the trust-region technique \cite{CGT2000,Deuflhard2004,Higham1999,LLT2007,Luo2010}. 
For this time-stepping scheme, it needs to select suitable a merit function and
construct an approximation model of the merit function. Here, we adopt the
residual $\|F(x)\|$ as the merit function and adopt
$\|F(x_{k}) + J(x_{k})s_{k}\|$ as the approximation model of
$\|F(x_{k} + s_{k})\|$. Thus, according to the following ratio:
\begin{align}
  \rho_k = \frac{\|F(x_{k})\|-\|F(x_{k} + s_{k})\|}
   {\|F(x_{k})\| - \|F(x_{k})+J(x_{k})s_{k})\|}, \label{RHOK}
\end{align}
we enlarge or reduce the time step $\Delta t_k$ at every iteration.
A particular adjustment strategy is given as follows:
\begin{align}
   \Delta t_{k+1} =
     \begin{cases}
    \gamma_1 \Delta t_k, &{\text{if} \;  \left|1- \rho_k \right| \le \eta_1,}\\
    \Delta t_k, &{\text{else if} \; \eta_1 < \left|1 - \rho_k \right| < \eta_2,}\\
    \gamma_2 \Delta t_k, &{\text{others},}
    \end{cases} \label{TSK1}
\end{align}
where the constants are selected as $\gamma_{1} = 2, \; \gamma_{2} = 0.5, \;
\eta_{1} = 0.25, \; \eta_{2} = 0.75$ according to our numerical experiments.

\vskip 2mm 

\begin{remark}
This new time-stepping scheme based on the trust-region updating
strategy has some advantages compared to the traditional line search strategy
\cite{Luo2005}. If we use the line search strategy and the damped Newton method
\eqref{NEWTON} to track the trajectory $z(t)$ of the continuous
Newton flow \eqref{DAEFLOW}, in order to achieve the fast convergence rate in
the steady-state phase, the time step size $\alpha_{k}$ of the damped Newton
method is tried from 1 and reduced by half with many times at every iteration.
Since the linear model $F(x_{k}) + J(x_{k})s_{k}$ may not approximate
$F(x_{k}+s_{k})$ well in the transient-state phase, the time step size $\alpha_{k}$
will be small. Consequently, the line search strategy consumes the unnecessary
trial steps in the transient-state phase. However, the selection scheme of the
time step size based on the trust-region updating strategy \eqref{RHOK}-\eqref{TSK1}
can overcome this shortcoming.
\end{remark}

\vskip 2mm 

According to the above discussions, we give the detailed implementation of
the regularization continuation Newton method with the residual trust-region 
time-stepping scheme for nonlinear equations in Algorithm \ref{CNMTR}.

\vskip 2mm 

\begin{algorithm}
   \renewcommand{\algorithmicrequire}{\textbf{Input:}}
   \renewcommand{\algorithmicensure}{\textbf{Output:}}
   \caption{Continuation Newton methods with the residual trust-region time-stepping scheme
   (The CNMTr method)}
   \label{CNMTR}
   \begin{algorithmic}[1]
      \REQUIRE ~~\\
        Function $F: \; \Re^{n} \to \Re^{n}$, the initial point $x_0$ \, (optional), 
        and the tolerance $\epsilon$ \, (optional).
	  \ENSURE ~~\\
        An approximation solution $x^{\ast}$ of nonlinear equations. \\
      \STATE Set the default $x_0 = \text{ones}(n, \; 1)$ and $\epsilon = 10^{-6}$ when 
      $x_{0}$ or $\epsilon$ is not provided by the calling subroutine. 

      \STATE Initialize the parameters: $\eta_{a} = 10^{-6}, \; \eta_1 = 0.25, \;
      \gamma_1 =2, \; \eta_2 = 0.75, \; \gamma_2 = 0.5, \; \text{maxit} = 400$. 
      \STATE Set $\Delta t_0 = 10^{-2}$, flag\_success\_trialstep = 1, $\text{itc} = 0, \; k = 0$.
      \STATE Evaluate $F_{k} = F(x_{k})$ and $J_{k} = J(x_{k})$. Compute the 
      residual $Res_0 = \|F(x_0)\|_{\infty}$.
      \WHILE{(itc $<$ maxit)}
         \IF{(flag\_success\_trialstep == 1)}
              \STATE Set itc = itc + 1.
              \STATE Compute $\text{Res}_{k} = \|F_{k}\|_{\infty}$.
              \IF{($\text{Res}_{k} < \epsilon$)}
                 \STATE break;
              \ENDIF
              \STATE Solve the linear system \eqref{RSMEDAE} to obtain the 
              Newton step $s_{k}^{P}$. 
         \ENDIF
         \STATE Compute $s_{k}  = {\Delta t_k}/{(1+\Delta t_k)} \, s_{k}^{P}$. 
         \STATE Set $x_{k+1} = x_k + s_k$.
         \STATE Evaluate $F(x_{k+1})$.            
         \IF {$\|F(x_k)\| < \|F(x_k) + J(x_k)s_k\|$}
           \STATE $\rho_{k} = -1$;
         \ELSE
           \STATE Compute the ratio $\rho_{k}$ from equation \eqref{RHOK}.
         \ENDIF
         \STATE Adjust the time-stepping size $\Delta t_{k+1}$ according to the
         trust-region updating strategy \eqref{TSK1}.
         \IF{$(\rho_{k} \ge \eta_{a})$}
               \STATE Accept the trial point $x_{k+1}$. Set flag\_success\_trialstep = 1.
           \ELSE
               \STATE Set $x_{k+1}  = x_{k}$,  $F_{k+1} = F_{k}$,  $s_{k+1}^{P} = s_{k}^{P}$, 
               flag\_success\_trialstep = 0.
           \ENDIF           
         \STATE Set $k \longleftarrow k+1$.
      \ENDWHILE
   \end{algorithmic}
\end{algorithm}

\vskip 2mm 

\section{Convergence analysis}
In this section, we discuss some theoretical properties of Algorithm \ref{CNMTR}.
Firstly, we estimate the lower bound of the predicted reduction
$\|F(x_k)\|-\|F(x_k) + J(x_k)s_{k}\|$, which is similar to the
estimate of the trust-region method for the unconstrained optimization problem
\cite{Powell1975}.

\begin{lemma} \label{LEMESTPRED}
Assume that it exists a positive constant $m$ such that
\begin{align}
  \|J(x_{k})y\| \ge m \|y\|, \; \forall y \in \Re^{n}, \; k = 0, \, 1, \,
  2, \, \ldots. \label{UBINVJ}
\end{align}
Furthermore, we suppose that $s_{k}$ is the solution of the regularization continuation
Newton method \eqref{RSMEDAE}-\eqref{XKSK}, where the regularization parameter
$\mu_{k}$ defined by equation \eqref{MUSW} and the constant $c_{\epsilon}$ satisfy
$\mu_{k} \le c_{\epsilon} < 0.5m$. Then, we have the following estimation
\begin{align}
   \|F(x_k)\| - \|F(x_k) + J(x_k)s_{k}\| \ge
    c_{r}\, \|F(x_k)\| \, {\Delta t_k}/{(1+\Delta t_k)},  \label{PRERED}
\end{align}
where the positive constant $c_{r}$ satisfies $0 < c_{r} < 1$.
\end{lemma}
\proof From equations \eqref{RSMEDAE}-\eqref{RSCORDAE}, we have
\begin{align}
  -J(x_k)s_{k} + \mu_{k} s_{k} =  ({\Delta t_k}/{(1+\Delta t_k)})F(x_k).
  \label{EQSK}
\end{align}
Thus, from equation \eqref{EQSK}, we obtain
\begin{align}
   & \|J(x_k)s_{k}  + F(x_k)\|   =
     \left\|\mu_{k} s_{k}+ {F(x_k)}/{(1+\Delta t_{k})}\right\| \nonumber \\
   & \quad = \left\|\mu_{k} ({\Delta t_k}/{(1+\Delta t_k)}) \, (-J(x_k)+\mu_{k}I)^{-1}F(x_k)
    + {F(x_k)}/{(1+\Delta t_{k})}\right\| \nonumber \\
   & \quad \le  ({1}/{(1+\Delta t_k)}) \, \left(\Delta t_k
    \left\|\left(- {J(x_k)}/{\mu_{k}}+ I\right)^{-1}\right\| +1\right)\|F(x_k)\|.
  \label{ESTRED}
\end{align}

\vskip 2mm

According to the definition of the induced matrix norm \cite{GV2013}, we have
\begin{align}
  & \left\|\left(-{J(x_k)}/{\mu_{k}} + I\right)^{-1}\right\|
  = \max_{z \neq 0} {\left\|\left(-{J(x_k)}/{\mu_{k}}
  + I\right)^{-1}z \right\|}/{\|z\|}  \nonumber \\
  & = \max_{ y\neq 0} \frac{\|y\|}{\left\|\left(-{J(x_k)}/{\mu_{k}}
  + I\right)y\right\|}
  = \frac{1}{\min_{ \|y\| = 1} \left\|\left(-{J(x_k)}/{\mu_{k}}
  + I\right)y\right\|}. \label{MATNORM}
\end{align}
On the other hand, when $\|y\| = 1$, from the nonsingular assumption
\eqref{UBINVJ} of matrix $J(x_{k})$, we have
\begin{align}
   & \left\|\left(- {J(x_k)}/{\mu_{k}} + I\right)y\right\|
    = \left\|- {J(x_k)y}/{\mu_{k}} + y \right\| \nonumber \\
   & \hskip 2mm  \ge {\|J(x_k)y\|}/{\mu_{k}} - \|y\|
    \ge  {m}/{\mu_{k}} - 1.    \label{LBMAT}
\end{align}
Thus, from the assumption $\mu_{k} \le c_{\epsilon} < 0.5m$ and equations
\eqref{MATNORM}-\eqref{LBMAT}, we have
\begin{align}
  \left\|\left(-{J(x_k)}/{\mu_{k}} + I\right)^{-1}\right\|
  \le {\mu_{k}}/{(m - \mu_{k})}
  \le {c_{\epsilon}}/{(m - c_{\epsilon})}. \label{LBMATN}
\end{align}

\vskip 2mm

By substituting inequality \eqref{LBMATN} into inequality \eqref{ESTRED}, we have
\begin{align}
   \|J(x_k)s_{k}  + F(x_k)\| \le (\left(1 + {c_{\epsilon}\Delta t_{k}}
   /{(m - c_{\epsilon})}\right)/{(1+\Delta t_{k})})\|F(x_k)\|.
    \nonumber
\end{align}
That is to say, we obtain 
\begin{align}
  \|F(x_k)\| - \|F(x_k)+J(x_k)s_{k}\|
  \ge \left(\frac{m - 2 c_{\epsilon}}{m - c_{\epsilon}}\right)
  \frac{\Delta t_{k}} {1+\Delta t_{k}}\|F(x_k)\|. \label{PREDEST}
\end{align}
We set $c_{r} = (m - 2c_{\epsilon})/(m - c_{\epsilon})$ in the above inequality
\eqref{PREDEST}. Then, we obtain the estimation \eqref{PRERED}. \qed

\vskip 2mm

In order to prove that the sequence $\{\|F(x_k)\|\}$ converges to zero when $k$
tends to infinity, we also need to estimate the lower bound of the time step
size $\Delta t_{k}$.

\vskip 2mm

\begin{lemma} \label{DTBOUND}
Assume that $F: \; \Re^{n} \to \Re^{n}$ is continuously differentiable
and its Jacobian function  $J$ is Lipschitz continuous. That to say, 
it exists a positive number $L$ such that
\begin{align}
  \|J(x)-J(y)\| \le L\|x-y\|   \label{LIPCON}
\end{align}
holds for all $x, \, y \in \Re^{n}$. 
Furthermore, we suppose that the sequence $\{x_k\}$ is generated by Algorithm
\ref{CNMTR} and the nonsingular condition \eqref{UBINVJ} of matrix $J(x_{k})$
holds. Then, when the regularization parameter $\mu_{k}$ defined by equation
\eqref{MUSW} and the constant $c_{\epsilon}$ satisfy $\mu_{k} \le c_{\epsilon} < 0.5m$,
it exists a positive number $\delta_{\Delta t}$ such that
\begin{align}
  \Delta t_{k} \ge \delta_{\Delta t} > 0, \; k = 0, \, 1, \; 2, \dots, \label{DTGEPN}
\end{align}
where $\Delta t_{k}$ is adaptively adjusted by formulas \eqref{RHOK}-\eqref{TSK1}.
\end{lemma}

\proof From the Lipschitz continuous assumption \eqref{LIPCON} of $J(\cdot)$,
we have
\begin{align}
    & \left\|F(x_{k+1}) - F(x_k) - J(x_k)s_{k}\right\|
    = \left\|\int_{0}^{1}J(x_{k}+ ts_{k})s_{k}dt - J(x_k)s_{k}\right\|
    \nonumber \\
    & = \left\|\int_{0}^{1}(J(x_{k}+ ts_{k})-J(x_k))s_{k}dt \right\|
    \le \int_{0}^{1}\|(J(x_{k}+ ts_{k})-J(x_k))s_{k}\|dt \nonumber \\
    & \le \int_{0}^{1}\|J(x_{k}+ ts_{k})-J(x_k)\| \|s_{k}\|dt
    \le \int_{0}^{1}L \|s_{k}\|^{2} tdt = 0.5 \, L \|s_{k}\|^{2}.
    \label{LIPAPUB}
\end{align}
On the other hand, from equations \eqref{RSMEDAE}-\eqref{RSCORDAE}, we have
\begin{align}
   \|s_k\| & = ({\Delta t_k}/{(1+ \Delta t_k)})
   \left\|\left(-J(x_k)+\mu_{k}I \right)^{-1}F(x_k) \right\| \nonumber \\
   &  \le ({\Delta t_k}/{(1+ \Delta t_k)})
   \left\|\left(-J(x_k)+\mu_{k}I \right)^{-1}\right\| \|F(x_k)\|. \label{ESTSK}
\end{align}

\vskip 2mm 

Similarly to the estimation \eqref{LBMATN}, from the assumption
$\mu_{k} \le c_{\epsilon} < 0.5m$ and the nonsingular assumption \eqref{UBINVJ} of
$J(x_k)$, we have
\begin{align}
   \left\|(-J(x_k) + \mu_{k}I)^{-1}\right\| \le {1}/{(m-\mu_{k})}
   \le {1}/{(m - c_{\epsilon})}. \label{INJACMUUB}
\end{align}
Thus, from inequalities \eqref{LIPAPUB}-\eqref{INJACMUUB}, we obtain
\begin{align}
   \left\|F(x_{k+1}) - F(x_k) - J(x_k)s_{k}\right\|
   \le \frac{L}{2(m- c_{\epsilon})^{2}}
   \left(\frac{\Delta t_k}{1+\Delta t_k}\right)^{2}\|F(x_k)\|^{2}.
  \label{LINAPUPB}
\end{align}

\vskip 2mm 

From the definition \eqref{RHOK} of $\rho_{k}$, the estimation
\eqref{PREDEST}, and inequality \eqref{LINAPUPB}, we have
\begin{align}
  & |\rho_{k} - 1| = \left|\frac{\|F(x_k)\|-\|F(x_{k+1})\|}
  {\|F(x_k)\|-\|F(x_k)+J(x_k)s_{k}\|}-1 \right| \nonumber \\
    & \hskip 2mm \le \frac{\|F(x_{k+1}) - F(x_k) - J(x_k)s_{k}\|}
  {\|F(x_k)\|-\|F(x_k)+J(x_k)s_{k}\|}
   \le \frac{L}{2(m-2c_{\epsilon})^{2}}
   \left(\frac{\Delta t_k}{1+\Delta t_k}\right) \|F(x_k)\|.
   \label{ESTRHO}
\end{align}
According to Algorithm \ref{CNMTR}, we know that the sequence $\{\|F(x_k)\|\}$
is monotonically decreasing. Consequently, we have $\|F(x_k)\| \le \|F(x_0)\|, \;
k = 1, \, 2, \, \dots$. We set
\begin{align}
    \bar{\delta}_{\Delta t} \triangleq {2(m-2c_{\epsilon})^{2}\eta_{1}}/{(\|F(x_0)\|L)} .
  \label{PARDELAT}
\end{align}

\vskip 2mm 

Assume that $K$ is the first index such that $\Delta t_{K} \le \bar{\delta}_{\Delta t}$.
Then, from inequalities \eqref{ESTRHO}-\eqref{PARDELAT}, we obtain
$|\rho_{K} - 1| < \eta_1$. Consequently, $\Delta t_{K+1}$ will be greater than
$\Delta t_{K}$ according to the adaptive adjustment scheme \eqref{TSK1}.
We set $\delta_{\Delta t} = \min \{\Delta t_K, \bar{\delta}_{\Delta t} \}$. Then,
$\Delta t_k \ge \delta_{\Delta t}$ holds for all $k = 0, \, 1, \, 2, \, \dots$. \qed

\vskip 2mm 

By using the estimation results of Lemma \ref{LEMESTPRED} and Lemma \ref{DTBOUND}, we can prove
that the sequence $\{\|F(x_k)\|\}$ converges to zero when $k$ tends to infinity.

\vskip 2mm 

\begin{theorem} \label{THECOVFXK}
Assume that $F: \; \Re^{n} \to \Re^{n}$ is continuously differentiable
and its Jacobian function $J$ satisfies the Lipschitz condition \eqref{LIPCON}.
Furthermore, we suppose that the sequence $\{x_k\}$ is generated by Algorithm
\ref{CNMTR} and $J(x_{k})$ satisfies the nonsingular assumption \eqref{UBINVJ}.
Then, when the regularization parameter $\mu_{k}$ defined by equation \eqref{MUSW}
and the constant $c_{\epsilon}$ satisfy $\mu_{k} \le c_{\epsilon} < 0.5m$, we have
\begin{align}
    \lim_{k \to \infty} \inf \; \|F(x_k)\| = 0. \label{FKTOZ}
\end{align}
\end{theorem}

\vskip 2mm 

\proof According to Algorithm \ref{CNMTR} and inequality \eqref{ESTRHO}, we
know that there exists an infinite subsequence $\{x_{k_{l}}\}$ such that
\begin{align}
   \frac{\|F(x_{k_{l}})\| - \|F(x_{k_l}+s_{k_l})\|}{\|F(x_{k_l})\|-
    \|F(x_{k_l})+J(x_{k_l})s_{k_l}\|} \ge \eta_{a}, \;
    l = 1, \, 2, \, \ldots . \label{ASUBSQ}
\end{align}
Otherwise, all steps are rejected after a given iteration index, then the time step 
size $\Delta t_{k}$ will keep decreasing, which contradicts equation \eqref{DTGEPN}.

\vskip 2mm

From inequalities \eqref{PRERED}, \eqref{DTGEPN} and \eqref{ASUBSQ}, we have
\begin{align}
    \|F(x_{k_{l}})\| - \|F(x_{k_l} + s_{k_l})\| \ge \frac{\eta_{a}
     c_{r}\Delta t_{k_l}}{(1+\Delta t_{k_l})}\|F(x_{k_l})\|
     \ge \frac{\eta_{a} c_{r}\delta_{\Delta t}}{(1+\delta_{\Delta t})}\|F(x_{k_l})\|. \label{FK1GE}
\end{align}
Therefore, from equation \eqref{FK1GE} and $\|F(x_{k+1})\| \le \|F(x_k)\|$, we
have
\begin{align}
    & \|F(x_0)\| \ge \|F(x_0)\| - \lim_{k \to \infty} \|F(x_k)\| =
    \sum_{k = 0}^{\infty} (\|F(x_k)\| - \|F(x_{k+1})\|) \nonumber \\
    & \hskip 2mm \ge \sum_{l = 0}^{\infty} (\|F(x_{k_l})\| - \|F(x_{k_l}+s_{k_l})\|)
    \ge \frac{\eta_{a}c_{r}\delta_{\Delta t}}{1+ \delta_{\Delta t}} \sum_{l=0}^{\infty}\|F(x_{k_l})\|. 
   \label{ESTSUBSQ}
\end{align}
Consequently, from inequality \eqref{ESTSUBSQ}, we obtain
\begin{align}
    \lim_{k_l \to \infty} \|F(x_{k_l})\| = 0. \label{SUBFXCON}
\end{align}
That is to say, the result \eqref{FKTOZ} is true.  Furthermore, from 
$\|F(x_{k+1})\| \le \|F(x_{k})\|$  and equation \eqref{SUBFXCON}, 
it is not difficult to know $\lim_{k \to \infty} \|F(x_{k})\| = 0$. \qed

\vskip 2mm

Under the nonsingular assumption of $J(x^{\ast})$ and the local Lipschitz continuity
\eqref{LIPCON} of $J(\cdot)$, we analyze the local superlinear convergence rate 
of Algorithm \ref{CNMTR} near the solution $x^{\ast}$ as follows. For convenience, 
we define the neighbourhood $B_{\delta}(x^{\ast})$ of $x^{\ast}$ as 
\begin{align}
B_{\delta}(x^{\ast})= \{x: \|x - x^{\ast}\|\le \delta \}. \nonumber 
\end{align}

\vskip 2mm

\begin{theorem} \label{LOCSCON}
Assume that $F: \; \Re^{n} \to \Re^{n}$ is continuously differentiable and
$F(x^{\ast}) = 0$. Furthermore, we suppose that $J$ satisfies the local
Lipschitz continuity \eqref{LIPCON} around $x^{\ast}$ and the nonsingular condition
\eqref{UBINVJ} when $x \in B_{\delta}(x^{\ast})$. Then, when the regularization
parameter $\mu_{k}$ defined by equation \eqref{MUSW} and the constant
$c_{\epsilon}$ satisfy $\mu_{k} \le c_{\epsilon} < 0.5m$, there exists a neighborhood
$B_{r}(x^{\ast})$ such that the sequence $\{x_k\}$ generated by
Algorithm \ref{CNMTR} with $x_{0} \in B_{r}(x^{\ast})$ superlinearly converges
to $x^{\ast}$.
\end{theorem}
\proof The framework of its proof can be roughly described as follows.
Firstly, we prove that the sequence $\{x_k\}$ linearly converges to
$x^{\ast}$ when $x_{0}$ gets close enough to $x^{\ast}$. Then, we prove
$\lim_{k \to \infty} \Delta t_k = \infty$. Finally, we prove that
the search step $s_{k}$ approximates the Newton step $s_{k}^{N}$. Consequently,
the sequence $\{x_{k}\}$ superlinearly converges to $x^{\ast}$.

\vskip 2mm

Firstly, similarly to the estimation \eqref{LBMATN}, from the assumption
$\mu_{k} \le c_{\epsilon} < 0.5m$, we obtain
\begin{align}
   \left\|(\mu_{k}I - J(x_{k}))^{-1}\right\| \le {1}/{(m - c_{\epsilon})},
   \; \forall x_{k} \in B_{\delta}(x^{\ast}), \; k = 0, \, 1, \, 2, \ldots.
   \label{INJACUPB}
\end{align}
We denote $e_{k} = x_{k} - x^{\ast}$. From equations \eqref{RSMEDAE}-\eqref{RSCORDAE},
we have
\begin{align}
   e_{k+1} & = e_{k} + s_{k} = e_{k} + \frac{\Delta t_{k}}{1+ \Delta t_{k}}
   \left(\mu_{k}I -J(x_k)\right)^{-1} (F(x_k) - F(x^{\ast}))
   \nonumber \\
   & = e_{k} + \frac{\Delta t_{k}}{1+ \Delta t_{k}}
   \left(\mu_{k}I -J(x_k)\right)^{-1}\int_{0}^{1} J(x^{\ast} + t e_{k})e_{k}dt.
   \label{EKSK}
\end{align}
By rearranging the above equation \eqref{EKSK}, we obtain
\begin{align}
   e_{k+1} = \frac{1}{1+\Delta t_k}e_{k} + \frac{\Delta t_{k}}{1+ \Delta t_{k}}
  \left(\mu_{k}I -J(x_k)\right)^{-1}\int_{0}^{1} \left(J(x^{\ast} + t e_{k})
  -J(x_k) + \mu_{k}I \right)e_{k}dt.  \nonumber
\end{align}
By using the Lipschitz continuous assumption \eqref{LIPCON} of $J(\cdot)$, the estimation
\eqref{INJACUPB}, and the assumption $\mu_{k} \le c_{\epsilon} < 0.5m$, we have
\begin{align}
   & \|e_{k+1}\|  \le {\|e_k\|}/{(1+\Delta t_k)} \nonumber \\
     & \hskip 2mm + ({\Delta t_k}/{(1+ \Delta t_k)})
     \left\| \left(\mu_{k}I -J(x_k)\right)^{-1}\right\| \int_{0}^{1}
     \left(\left\|J(x^{\ast} + t e_{k}) -J(x_k)\right\| + \mu_{k} \right)\|e_{k}\|dt
     \nonumber \\
   & \hskip 2mm  \le {\|e_k\|}/{(1+\Delta t_k)} + ({\Delta t_k}/{((1+ \Delta t_k)(m - \mu_{k}))})
    \left(\mu_{k} + 0.5 L\|e_k\|\right)\|e_k\|
    \nonumber \\
   & \hskip 2mm  = \frac{1+\frac{1}{m- \mu_{k}}\left(\mu_{k} + 0.5L\|e_k\|\right){\Delta t_k}}
     {1+\Delta t_k}\|e_k\|
     \le \frac{1+\frac{1}{m- c_{\epsilon}}\left(c_{\epsilon} + 0.5L\|e_k\|\right){\Delta t_k}}
     {1+\Delta t_k}\|e_k\|. \label{EK1EKNS}
\end{align}

\vskip 2mm

We denote
\begin{align}
   q_{k} \triangleq \frac{1+ 
   \left(c_{\epsilon} + 0.5L\|e_k\|\right){\Delta t_k}/({m - c_{\epsilon}})}
     {1+\Delta t_k}, \label{QK}
\end{align}
and select $x_{0} \in B_{\delta}(x^{\ast})$  to satisfy
\begin{align}
  \|e_{0}\| < {(m - 2 c_{\epsilon})}/{L}.   \label{EKLESS1}
\end{align}
We set $r = \min \{\delta, (m - 2c_{\epsilon})/L\}$.  When $x_{0} \in B_{r}(x^{\ast})$,
from equations \eqref{EK1EKNS}-\eqref{EKLESS1} and the assumption $c_{\epsilon} < 0.5m$,
by induction, we have
\begin{align}
   \|e_{k+1}\| \le q_{k} \|e_k\|, \;
    q_{k} < \frac{1 + {0.5 \Delta t_k \, m}/{(m - c_{\epsilon})}}{1+\Delta t_k} < 1,
   \; k = 0, \, 1, \, \ldots.   \label{CONEK1}
\end{align}

\vskip 2mm

It is not difficult to know that $f(t) \triangleq (1+\alpha t)/(1+t)$ is monotonically
decreasing when $0 \le \alpha < 1$. Thus, from the estimation \eqref{DTGEPN} of
the time step size $\Delta t_{k}$ and inequality \eqref{CONEK1}, we obtain
\begin{align}
   \|e_{k+1}\| \le q_{k} \|e_{k}\| \le q \|e_{k}\|, \;
   q \triangleq \frac{1 + {0.5 \delta_{\Delta t} m}/{(m - c_{\epsilon})}}
   {1+\delta_{\Delta t}} < 1.     \nonumber
\end{align}
Therefore, we have
\begin{align}
   \|e_{k+1}\| \le q^{k} \|e_{0}\| \to 0, \; \text{when}   \; k \to \infty.
   \label{EKTOZOR}
\end{align}
That is to say, we obtain $\lim_{k \to \infty} x_k = x^{\ast}$.

\vskip 2mm

Secondly, from equations \eqref{RSMEDAE}-\eqref{RSCORDAE} and inequality
\eqref{INJACUPB}, we have
\begin{align}
  \|s_k\| & =  ({\Delta t_k}/{(1+\Delta t_k)}) \left\|(-J(x_{k}) + \mu_{k}I)^{-1}F(x_k)\right\| 
   \nonumber \\
  & \le \frac{\Delta t_k}{1+\Delta t_k} \left\|(-J(x_{k}) + \mu_{k}I)^{-1}\right\| \|F(x_k)\|
  \le \frac{1}{m - c_{\epsilon}}\frac{\Delta t_k}{1+\Delta t_k}\|F(x_k)\|. \label{SKUPB}
\end{align}
Similarly to the estimation \eqref{ESTRHO}, from the definition
\eqref{RHOK} of $\rho_{k}$, inequalities \eqref{PREDEST} and \eqref{SKUPB}, we have
\begin{align}
  |\rho_{k} - 1| & = \left|\frac{\|F(x_k)\|-\|F(x_{k+1})\|}
  {\|F(x_k)\|-\|F(x_k)+J(x_k)s_{k}\|}-1 \right| \nonumber \\
   & \le \frac{0.5L}{(m - 2c_{\epsilon})^{2}}
   \left(\frac{\Delta t_k}{1+\Delta t_k}\right) \|F(x_k)\|
   \le \frac{0.5L}{(m - 2c_{\epsilon})^{2}} \|F(x_k)\|.  \label{ESTRHOEN}
\end{align}
Since $\{\|F(x_k)\|\}$ is monotonically decreasing and $\lim_{k \to \infty} x_{k} = x^{\ast}, \; 
F(x^{\ast}) = 0$, we can select a sufficiently large number $K$ such that
\begin{align}
   \|F(x_{k})\| \le \frac{2\eta_{1}(m - 2c_{\epsilon})^{2}}{L}, \; \text{when}
   \;   k \ge K.   \label{FXKLES}
\end{align}
From inequalities \eqref{ESTRHOEN}-\eqref{FXKLES}, we have
\begin{align}
   |\rho_{k} - 1| \le \eta_1, \; \text{when}   \; k \ge K.  \nonumber
\end{align}
This means $\Delta t_{k+1} = \gamma_{1} \Delta t_{k}$ when $k \ge K$,
according to the time-stepping scheme \eqref{TSK1}. That is to say, we have
\begin{align}
   \lim_{k \to \infty} \Delta t_k = \infty. \label{DELTOINF}
\end{align}

\vskip 2mm

Finally, since $\lim_{k \to \infty} \Delta t_k = \infty$, we can select a
sufficiently large number $K_{\mu}$ such that $1/\Delta t_k < c_{\epsilon}$
when $k \ge K_{\mu}$. Consequently, from the definition \eqref{MUSW} of
the regularization parameter $\mu_{k}$, we obtain $\mu_{k} = 1/\Delta t_{k}$
when $k \ge K_{\mu}$. By substituting  it into inequality \eqref{EK1EKNS},
we have
\begin{align}
    & \frac{\|e_{k+1}\|}{\|e_{k}\|}   \le
    \frac{1}{1+\Delta t_k}  + \frac{\Delta t_k}{1+ \Delta t_k}
     \frac{1}{m - \mu_{k}}\left(\mu_{k} + 0.5L\|e_k\|\right)
     \nonumber \\
     & \hskip 2mm  =  \frac{1}{1+\Delta t_k} + \frac{\Delta t_k}{1+\Delta t_k}
     \frac{1}{m-{1}/{\Delta t_k}}
    \left(1/{\Delta t_k} + 0.5 L \|e_k\|\right), \; \text{when} \; k \ge K_{\mu}.
    \label{EK1VSEK}
\end{align}
From equations \eqref{EKTOZOR} and \eqref{DELTOINF}, we know
$\lim_{k \to \infty} \|e_{k}\| = 0$ and $\lim_{k \to \infty} \Delta t_{k} = \infty$,
respectively. Therefore,  by combining them with inequality \eqref{EK1VSEK}, we
obtain
\begin{align}
   \lim_{k \to \infty} \frac{\|e_{k+1}\|}{\|e_{k}\|} = 0. \nonumber
\end{align}
That is to say, the sequence $\{x_{k}\}$ superlinearly converges to $x^{\ast}$.
\qed

\vskip 2mm

For the real-world problem, the singularity of $J(x)$ may arise from
the linear conservation law such as the conservation of mass or
the conservation of charge \cite{Luo2009,Shampine1998,Shampine1999,SGT2003}.
In the rest of this section, we analyze convergence properties of Algorithm
\ref{CNMTR} when $J(x)$ is singular. Similarly to the standard assumption of
the nonlinear dynamical system, we suppose that $J(\cdot)$ satisfies the one-sided
Lipschitz condition (see p. 303, \cite{Deuflhard2004} or p. 180, \cite{HW1996})
as follows:
\begin{align}
   y^{T}J(x)y \le - \nu \|y\|^2, \; \text{for} \; y \in {S}_{c} =
   \left \{y|c^{T}y = 0 \right \}, \; \nu > 0, \label{OSLC}
\end{align}
where the constant vector $c$ satisfies $c^{T}F(x) = 0, \; \forall x\in \Re^{n}$.
The positive number $\nu$ is called the one-sided Lipschitz constant. Under the assumption
of the one-sided Lipschitz condition \eqref{OSLC}, we know that matrix $(\mu I - J(x))$ is
nonsingular when $\mu > 0$. We state it as the following property \ref{PRONONSIG}.

\vskip 2mm

\begin{property} \label{PRONONSIG}
Assume that $J(\cdot)$ satisfies the one-sided Lipschitz condition \eqref{OSLC}. Then, matrix $(\mu I - J(x))$
is nonsingular when $\mu >0$, and the solution $s_{k}$ of equations \eqref{RSMEDAE}-\eqref{RSCORDAE} satisfies
$c^{T}s_{k} = 0$.
\end{property}
\proof We prove it by contradiction. If we assume that matrix $(\mu I - J(x))$ is
singular, there exists a nonzero vector $y$ such that
\begin{align}
     (\mu I - J(x)) y = 0. \label{ASINJX}
\end{align}
Consequently, from the assumption $c^{T}F(x) = 0$, we have
\begin{align}
    c^{T}y =  c^{T}(J(x)y)/\mu = \left(c^{T}J(x)\right)y/\mu = 0. \nonumber
\end{align}
Thus, from the one-sided Lipschitz condition \eqref{OSLC} and $\mu > 0, \; \nu > 0$,
we obtain
\begin{align}
  y^{T}(\mu I - J(x))y
  = \mu \|y\|^{2} - y^{T}J(x)y \ge (\mu + \nu)\|y\|^{2} > 0, \nonumber
\end{align}
which contradicts the assumption \eqref{ASINJX}. Therefore, matrix $(\mu I - J(x))$
is nonsingular.

\vskip 2mm

From equations \eqref{RSMEDAE}-\eqref{RSCORDAE}, we have
\begin{align}
   \left(\mu_{k} I - J(x_k) \right)s_{k}
   = ({\Delta t_k}/{(1+\Delta t_k)})F(x_k) .   \label{SK}
\end{align}
By combining it with the assumption $c^{T}F(x_k) = 0$, we obtain
\begin{align}
   c^{T}(\mu_{k}I - J(x_k))s_{k} = ({\Delta t_k}/{(1+\Delta t_k)}) c^{T}F(x_k)
    = 0. \label{CSK}
\end{align}
Therefore, by substituting $c^{T}J(x_k)= 0$ into equation \eqref{CSK}, we have
\begin{align}
   \mu_{k} c^{T}s_{k} = c^{T}J(x_k)s_{k} = \left(c^{T}J(x_k)\right)s_{k}
   = 0. \label{SKNULLC}
\end{align}
That is to say, we obtain $c^{T}s_k = 0$.  \qed

\vskip 2mm

Similarly to the estimation \eqref{PRERED} of $\|F(x_k)\| - \|F(x_k)+J(x_k)s_k\|$
for the nonsingular Jacobian matrix $J(x_{k})$, we also have its lower-bounded
estimation when $J(x_{k})$ is singular, $k = 0, \, 1, \, 2, \, \ldots$.

\vskip 2mm

\begin{lemma} \label{LEMSJEPRED}
Assume that $J(x_{k})$ satisfies the one-sided Lipschitz condition \eqref{OSLC}
and $s_{k}$ is the solution of equations \eqref{RSMEDAE}-\eqref{XKSK}. Then, we have
\begin{align}
   \|F(x_k)\| - \|F(x_k) + J(x_k)s_{k}\| \ge
   c_{s}({\Delta t_k}/{(1+\Delta t_k)})\|F(x_k)\|,  \label{SJPRERED}
\end{align}
where the positive constant $c_{s}$ satisfies $0 < c_{s} < 1$.
\end{lemma}

\vskip 2mm

\proof From Property \ref{PRONONSIG}, we know that matrix $(\mu I - J(x_k))$
is nonsingular and $s_k$ satisfies $c^{T}s_{k} = 0$. From equations
\eqref{RSMEDAE}-\eqref{RSCORDAE} and the Cauchy-Schwartz inequality
$|x^{T}y| \le \|x\| \|y\|$, we have
\begin{align}
   \mu_{k} \|s_k\|^{2} - s_{k}^{T}J(x_k)s_{k} = {s_{k}^{T}F(x_k) \Delta t_k}/{(1+\Delta t_k)}
     \le ({\Delta t_k}/{(1+\Delta t_k)})\|s_k\| \|F(x_k)\|.
   \label{SKJKLEFK}
\end{align}
By substituting one-sided Lipschitz condition \eqref{OSLC} into
equation \eqref{SKJKLEFK}, we obtain
\begin{align}
   (\mu_{k} + \nu) \|s_k\|^{2} \le \mu_{k} \|s_k\|^{2} - s_{k}^{T}J(x_k)s_{k}
   \le ({\Delta t_k}/{(1+\Delta t_k)})\|s_k\| \|F(x_k)\|. \nonumber
\end{align}
Consequently, we have
\begin{align}
   \|s_k\| \le \frac{1}{\mu_{k} + \nu}\frac{\Delta t_k}{1+\Delta t_k} \|F(x_k)\|.
   \label{SKLESFK}
\end{align}

From equations \eqref{RSMEDAE}-\eqref{RSCORDAE} and \eqref{SKLESFK}, we
have
\begin{align}
   & \|F(x_k) + J(x_k)s_{k}\|  = \left\|\mu_{k} s_k + \frac{1}{1+\Delta t_k}F(x_k)\right\|
   \le \mu_{k} \|s_k\| + \frac{1}{1+\Delta t_k}\|F(x_k)\|
   \nonumber \\
   & \hskip 2mm \le \frac{\mu_{k}}{\mu_{k} + \nu} \frac{\Delta t_k}{1+\Delta t_k}\|F(x_k)\|
   + \frac{1}{1+\Delta t_k}\|F(x_k)\|. \label{LINAPPFK1}
\end{align}
From the definition \eqref{MUSW} of the parameter $\mu_{k}$, we know
$\mu_{k} \le c_{\epsilon}$. By substituting it into inequality \eqref{LINAPPFK1},
we obtain
\begin{align}
    & \|F(x_k)\| - \|F(x_k) + J(x_k)s_{k}\|
    \ge \frac{\nu}{\mu_{k} + \nu} \frac{\Delta t_k}{1+\Delta t_k}\|F(x_k)\|
    \nonumber \\
   & \quad \ge \frac{\nu}{c_{\epsilon} + \nu} \frac{\Delta t_k}{1+\Delta t_k}\|F(x_k)\|.
   \label{FKMINFKJSK}
\end{align}
We set $c_{s} = \nu/(c_{\epsilon}+\nu)$. Then, from equation 
\eqref{FKMINFKJSK}, we obtain the estimation \eqref{SJPRERED}. \qed

\vskip 2mm

Similarly to the lower-bounded estimation \eqref{DTGEPN} of the time step size
$\Delta t_{k}$ for the nonsingular Jacobian matrix $J(x_{k})$, we also have its
lower-bounded estimation when $J(x_{k})$ is singular, $k = 0, \, 1, \, 2, \, \ldots$.

\vskip 2mm

\begin{lemma} \label{SDTBOUND}
Assume that $F: \; \Re^{n} \to \Re^{n}$ is continuously differentiable
and $J(\cdot)$ satisfies the Lipschitz continuity \eqref{LIPCON}
and the one-sided Lipschitz condition \eqref{OSLC}. The sequence $\{x_k\}$ is
generated by Algorithm \ref{CNMTR}. Then, there exists a positive
$\delta_{s}$ such that
\begin{align}
   \Delta t_{k} \ge \delta_{s}> 0, \; k = 1, \; 2, \ldots, \label{SDTGEPN}
\end{align}
where $\Delta t_{k}$ is adaptively updated by the trust-region updating strategy 
\eqref{RHOK}-\eqref{TSK1}.
\end{lemma}

\proof From the Lipschitz continuity \eqref{LIPCON} of $J(\cdot)$, we have
\begin{align}
   & \left\|F(x_{k+1}) - F(x_k) - J(x_k)s_{k}\right\|
   = \left\|\int_{0}^{1}(J(x_{k}+ ts_{k})-J(x_k))s_{k}dt \right\| \nonumber \\
   & \hskip 2mm  \le \int_{0}^{1}\|J(x_{k}+ ts_{k})-J(x_k)\| \|s_{k}\|dt
   \le \int_{0}^{1}L \|s_{k}\|^{2} tdt = 0.5 L \|s_{k}\|^{2}.
  \label{SLIPAPUB}
\end{align}
By substituting the estimation \eqref{SKLESFK} of $s_k$ into inequality
\eqref{SLIPAPUB}, we obtain
\begin{align}
   \left\|F(x_{k+1}) - F(x_k) - J(x_k)s_{k}\right\| \le
   \frac{0.5 L}{(\mu_{k} + \nu)^2}
   \left(\frac{\Delta t_k}{1+\Delta t_k}\right)^{2} \|F(x_k)\|^{2}.
   \label{FXK1JXKSK}
\end{align}
Thus, from the definition \eqref{RHOK} of $\rho_{k}$, inequalities
\eqref{FKMINFKJSK} and \eqref{FXK1JXKSK}, we obtain
\begin{align}
   & |\rho_{k} - 1| = \left|\frac{\|F(x_k)\|-\|F(x_{k+1})\|}
   {\|F(x_k)\|-\|F(x_k)+J(x_k)s_{k}\|}-1 \right|
    \le \frac{\|F(x_{k+1}) - F(x_k) - J(x_k)s_{k}\|}
   {\|F(x_k)\|-\|F(x_k)+J(x_k)s_{k}\|} \nonumber \\
   & \quad \le \frac{L}{2\nu(\mu_{k} + \nu)}
   \left(\frac{\Delta t_k}{1+\Delta t_k}\right) \|F(x_k)\|
   \le \frac{L}{2\nu^{2}} \left(\frac{\Delta t_k}{1+\Delta t_k}\right) \|F(x_k)\|.
   \label{SESTRHO}
\end{align}

According to Algorithm \ref{CNMTR}, we know that the sequence $\{\|F(x_k)\|\}$
is monotonically decreasing. Consequently, we have $\|F(x_k)\| \le \|F(x_0)\|$,
\; $k = 1, \, 2, \, \dots $. We set
\begin{align}
    \bar{\delta}_{s} = {2\nu^{2}\eta_{1}}/{(\|F(x_0)\|L)}.
   \label{SPARDELAT}
\end{align}
If we assume that $K$ is the first index such that $\Delta t_{K} \le \bar{\delta}_{s}$, 
then, from inequalities \eqref{SESTRHO}-\eqref{SPARDELAT}, we obtain
$|\rho_{K} - 1| < \eta_1$. Consequently, $\Delta t_{K+1}$ will be greater than
$\Delta t_{K}$ according to the time-stepping scheme \eqref{TSK1}. Set $\delta_{s} =
\min \{\Delta t_K, \bar{\delta}_{s} \}$. Then, we have $\Delta t_k \ge
\delta_{s}, \; k = 0, \, 1, \, 2, \, \dots$. \qed

\vskip 2mm

Now, from Lemma \ref{LEMSJEPRED} and Lemma \ref{SDTBOUND}, we know that the
sequence $\{\|F(x_k)\|\}$ converges to zero when $k$ tends to infinity and its
proof is similar to the proof of Theorem \ref{THECOVFXK}. We state it as the
following theorem \ref{STHECOVFXK} and omit its proof.

\vskip 2mm

\begin{theorem} \label{STHECOVFXK}
Assume that $F: \; \Re^{n} \to \Re^{n}$ is continuously differentiable and its Jacobian
function $J(\cdot)$ satisfies the Lipschitz continuity \eqref{LIPCON} and the one-sided Lipschitz
condition  \eqref{OSLC}. The sequence $\{x_k\}$ is generated by Algorithm
\ref{CNMTR}. Then, we have
\begin{align}
    \lim_{k \to \infty} \inf \; \|F(x_k)\| = 0. \label{SFKTOZ}
\end{align}
\end{theorem}

\begin{theorem} \label{SLOCSCON}
Assume that $F: \; \Re^{n} \to \Re^{n}$ is continuously differentiable
and its Jacobian function $J(\cdot)$ satisfies the Lipschitz continuity \eqref{LIPCON}
and the one-sided Lipschitz condition \eqref{OSLC}. Furthermore, we suppose
that the sequence $\{x_k\}$ is generated by Algorithm \ref{CNMTR} and its
subsequence $\{x_{k_i}\}$ converges to $x^{\ast}$. Then, the sequence $\{x_k\}$
superlinearly converges to $x^{\ast}$.
\end{theorem}

\proof The framework of its proof can be roughly described as follows.
Firstly, we prove $\lim_{k \to \infty} \Delta t_k = \infty$. Then, we prove that
the sequence $\{x_k\}$ linearly converges to $x^{\ast}$. Finally, we prove that
the search step $s_{k}$ approximates the Newton step $s_{k}^{N}$. Consequently,
the sequence $\{x_{k}\}$ superlinearly converges to $x^{\ast}$.

\vskip 2mm

From Property \ref{PRONONSIG}, we know that matrix $(\mu I - J(x_k))$
is nonsingular and $s_k$ satisfies $c^{T}s_{k} = 0$,  where the constant vector
$c$ satisfies $c^{T}F(x) = 0$ for all $x\in \Re^{n}$ and  $s_k$ is the solution
of equations \eqref{RSMEDAE}-\eqref{RSCORDAE}.

\vskip 2mm

Firstly, we prove that there exists an index $K$ such that $\Delta t_k$  will be
enlarged at every iteration when $k \ge K$. Consequently, we have
$\lim_{k \to \infty} \Delta t_k = \infty$. From the lower-bounded estimation
\eqref{FKMINFKJSK} of $F(x_{k}) - F(x_{k} + s_{k})$ and inequality \eqref{SESTRHO},
we have
\begin{align}
    & |\rho_{k} - 1| = \left|\frac{\|F(x_k)\|-\|F(x_{k+1})\|}
    {\|F(x_k)\|-\|F(x_k)+J(x_k)s_{k}\|}-1 \right| \nonumber \\
    & \quad \le \frac{L}{2\nu(\mu_{k}+\nu)}
   \left(\frac{\Delta t_k}{1+\Delta t_k}\right) \|F(x_k)\|
   \le \frac{L}{2\nu^2}\|F(x_k)\|.
   \label{SLSESTRHO}
\end{align}
Since the subsequence $\{x_{k_i}\}$ converges to $x^{\ast}$, there exists an index
$K_{F}$ such that
\begin{align}
    \|F(x_{K_F})\| \le {2\eta_1 \nu^2}/{L}.   \label{SFXKLES}
\end{align}
Furthermore, according to Algorithm \ref{CNMTR}, we know that the sequence
$\{\|F(x_k)\|\}$ is monotonically decreasing. Consequently, we have
$\|F(x_k)\| \le \|F(x_{K_F})\|$ when $k \ge K_F$. Thus, from inequalities
\eqref{SLSESTRHO}-\eqref{SFXKLES}, we have
\begin{align}
    |\rho_{k} - 1| \le \eta_1, \; \text{when}   \; k \ge K_F. \label{ROHLEETA1}
\end{align}
Consequently, according to the time-stepping scheme \eqref{TSK1}, we know that
$\Delta t_{k+1} = \gamma_{1} \Delta t_{k}$ when $k \ge K_F$. Therefore,
we obtain $\lim_{k \to \infty} \Delta t_k = \infty$.

\vskip 2mm

Secondly, we prove that the sequence $\{x_{k}\}$ linearly converges to $x^{\ast}$
as follows. We denote
\begin{align}
   e_{k} = x_{k} - x^{\ast}. \label{EK}
\end{align}
From equations \eqref{RSMEDAE}-\eqref{RSCORDAE} and \eqref{EK}, we have
\begin{align}
    e_{k+1} = e_k + s_k
    = e_k + ({\Delta t_k}/{(1+\Delta t_k)})\left(\mu_{k}I - J(x_k)\right)^{-1} F(x_k).
   \label{EK1EKSK}
\end{align}
By rearranging inequality \eqref{EK1EKSK}, we obtain
\begin{align}
     & \left(\mu_{k}I - J(x_k)\right)e_{k+1} = \left(\mu_{k}I - J(x_k)\right)e_{k}
     + ({\Delta t_k}/{(1+\Delta t_k)})(F(x_k) - F(x^{\ast}))  \nonumber \\
     & \quad = \mu_{k} e_k - {J(x_k)e_k}/{\Delta t_k}  
     + ({\Delta t_k}/{(1+\Delta t_k)})
     \int_{0}^{1}\left(J(x^{\ast}+te_{k})-J(x_k)\right)e_{k}dt. \label{EK1EK}
\end{align}

\vskip 2mm

Since $c^{T}s_{k}$ = 0 and the subsequence $\{x_{k_i}\}_{i=1}^{+\infty}$
converges to $x^{\ast}$, from equation \eqref{EK1EKSK}, we have
\begin{align}
   c^{T}e_{k+1} = c^{T}e_{k} + c^{T}s_{k} = c^{T}e_{k} =
   \cdots = c^{T}e_{k_i} \to 0, \; \text{when} \; i \to \infty.
 \nonumber
\end{align}
That is to say, we have $c^{T}e_{k} =0, \; k = 0, \, 1, \, \ldots$. Thus, from
the one-sided Lipschitz condition \eqref{OSLC} and the Cauchy-Schwartz inequality
$|x^{T}y| \le \|x\| \, \|y\|$, we have
\begin{align}
   &\|e_{k+1}\| \left\| \left(\mu_{k}I - J(x_k)\right)e_{k+1}\right\|
   \ge e_{k+1}^{T} \left(\mu_{k}I - J(x_k)\right)e_{k+1} \nonumber \\
   & = \mu_{k}e^{T}_{k+1}e_{k+1} - e^{T}_{k+1}J(x_k)e_{k+1} \ge
   (\mu_{k} + \nu)\|e_{k+1}\|^{2}. \label{EK1LEJK}
\end{align}
By rearranging inequality \eqref{EK1LEJK}, we obtain
\begin{align}
   \|e_{k+1}\| \le {\left\| \left(\mu_{k}I - J(x_k)\right)e_{k+1}\right\|}
   /{(\mu_{k} + \nu)}. \label{EK1LEUPB}
\end{align}

\vskip 2mm

From the continuity of $J$ at $x^{\ast}$, there exists the positive constants
$M$ and $\epsilon$ such that
\begin{align}
   \|J(x)\| \le M \; \text{when} \; \|x - x^{\ast}\| < \epsilon.
   \nonumber
\end{align}
Since the subsequence $\{x_{k_i}\}$ converges to $x^{\ast}$, there exists
$K_{1}$ such that $\|x_{K_1} - x^{\ast}\| < \epsilon$. By combining it with
$\lim_{k \to \infty} \Delta t_{k} = \infty$, we can select
a sufficiently large number $K_{2}$ such that $\Delta t_{K_2} \ge 4M/\nu$ and
$\|e_{K_2}\| \le 0.5\nu/L$. We set $K = \max \{K_1, \, K_2\}$.

\vskip 2mm 

From equation \eqref{EK1EK} and the Lipschitz continuity \eqref{LIPCON}, we have
\begin{align}
    &\|(\mu_{K}I - J(x_K))e_{K+1}\| \le
   \mu_{K} \|e_K\| + {\|J(x_K)\| \|e_K\|}/{\Delta t_K} 
   \nonumber \\
   & \hskip 2mm +  \int_{0}^{1}\|J(x^{\ast}+te_{K})-J(x_K)\| \|e_K\|dt \,
   {\Delta t_K}/{(1+\Delta t_K)}    \nonumber \\
   & \hskip 2mm \le \left(\mu_{K} + {M}/{\Delta t_K} \right) \|e_K\|
   +  \int_{0}^{1} L \|e_K\|^{2}t dt \le \left(\mu_{K} + {M}/{\Delta t_K} 
   + 0.5 L \|e_K\| \right)\|e_K\|.    \nonumber
\end{align}
By combining it with inequality \eqref{EK1LEUPB}, $\Delta t_{K} \ge (4M)/\nu$,
and $\|e_{K}\| \le \nu/(2L)$, we obtain
\begin{align}
  \|e_{K+1}\| \le \frac{\mu_{K} + {M}/{\Delta t_K}  + 0.5 L \|e_K\|}
   {\mu_{K} + \nu} \|e_K\|
    \le \frac{\mu_{K} + 0.5 \nu} {\mu_{K} + \nu} \|e_K\| < \|e_K\|< \epsilon.
    \label{SEK1LEK}
\end{align}
Therefore, by induction, we obtain
\begin{align}
     \|e_{k+1}\| \le \frac{\mu_{k} + {M}/{\Delta t_k}  + 0.5 L \|e_k\|}
   {\mu_{k} + \nu} \|e_k\| \le \frac{\mu_{k} + 0.5 \nu} {\mu_{k} + \nu} \|e_k\|,
   \; \text{when} \;  k \ge K.   \label{EK1LESUP}
\end{align}
Furthermore, from the definition \eqref{MUSW}, we know that $\mu_{k} < c_{\epsilon}$.
By substituting it into inequality \eqref{EK1LESUP}, we have
\begin{align}
    \|e_{k+1}\| \le q \|e_k\| \le \cdots \le q^{(k-K+1)} \|e_{K}\|, \; q \triangleq
    \frac{c_{\epsilon} + 1/2 \nu}{c_{\epsilon}+ \nu} < 1,
    \; \text{when} \; k \ge K.  \nonumber
\end{align}
Consequently, we obtain $\lim_{k \to \infty} \|e_k\| = 0$.

\vskip 2mm

Finally, we prove that the sequence $\{x_{k}\}$ superlinearly converges to
$x^{\ast}$. Since $\lim_{k \to \infty} \Delta t_k = \infty$, we can select a
sufficiently large number $K_{\mu}$ such that $1/\Delta t_k < c_{\epsilon}$ when
$k \ge K_{\mu}$. Thus, from the definition \eqref{MUSW} of the
regularization parameter $\mu_{k}$, we know $\mu_{k} = 1/\Delta t_{k}$ when
$k \ge K_{\mu}$. By substituting it into equation \eqref{EK1LESUP}, we obtain
\begin{align}
    \frac{\|e_{k+1}\|}{\|e_{k}\|} & \le
    \frac{\mu_{k} + {M}/{\Delta t_k}  + 0.5 L \|e_k\|}{\mu_{k} + \nu}
    = \frac{{1}/{\Delta t_k} + {M}/{\Delta t_k}  + 0.5 L \|e_k\|}
   {{1}/{\Delta t_k} + \nu}.    \label{SUPEK1EK}
\end{align}
Consequently, from $\lim_{k \to \infty} \Delta t_{k} = \infty$, 
$\lim_{k \to \infty} \|e_{k}\| = 0$, and equation \eqref{SUPEK1EK},  we have 
$\lim_{k \to \infty}{\|e_{k+1}\|}/{\|e_{k}\|} = 0$. 
That is to say, the sequence $\{x_{k}\}$ superlinearly converges to $x^{\ast}$. \qed

\vskip 2mm

\section{Numerical Experiments}

\vskip 2mm

In this section, for some real-world equilibrium problems and the classical test
problems of nonlinear equations, we test the performance of Algorithm \ref{CNMTR}
(CNMTr) and compare it with the trust-region method (the built-in subroutine
fsolve.m of the MATLAB environment \cite{MATLAB,More1978}) and the homotopy
methods (HOMPACK90 \cite{WSMMW1997}, and NAClab
\cite{LLT2008,ZL2013,Zeng2019}).

\vskip 2mm

HOMPACK90 \cite{WSMMW1997} is a classical homotopy method implemented by
fortran 90 for nonlinear equations and it is very popular in engineering fields.
Another state-of-the-art homotopy method is the built-in subroutine psolve.m of
the NAClab environment \cite{LLT2008,ZL2013}. Since psolve.m only solves the
polynomial systems, we replace psolve.m with its subroutine GaussNewton.m
(the Gauss-Newton method) for non-polynomial systems. Therefore, we compare these
two homotopy methods with Algorithm \ref{CNMTR}, too.

\vskip 2mm

We collect $26$ test problems of nonlinear equations, some of which come
from the equilibrium problems of chemical reactions \cite{HW1996,KLH1982,Robertson1966,VL1994},
and some of which come from the classical test problems
\cite{Deuflhard2004,DS2009,LUL1994,MGH1981,NW1999}. Their simple descriptions
are given by Table \ref{TABPROB}. Their dimensions vary from $1$ to $3000$. The Jacobian
matrix $J(\cdot)$ of $F(\cdot)$ is singular for some test problems. The codes are executed
by a HP Pavilion notebook with an Intel quad-core CPU. The termination condition
is given by
\begin{align}
   \|F(x^{it})\|_{\infty} \le 10^{-12}. \label{TERMCON}
\end{align}

\vskip 2mm

\begin{figure}
   \centering
   \begin{minipage}{.5\textwidth}
       \centering
       \includegraphics[width=.9\linewidth]{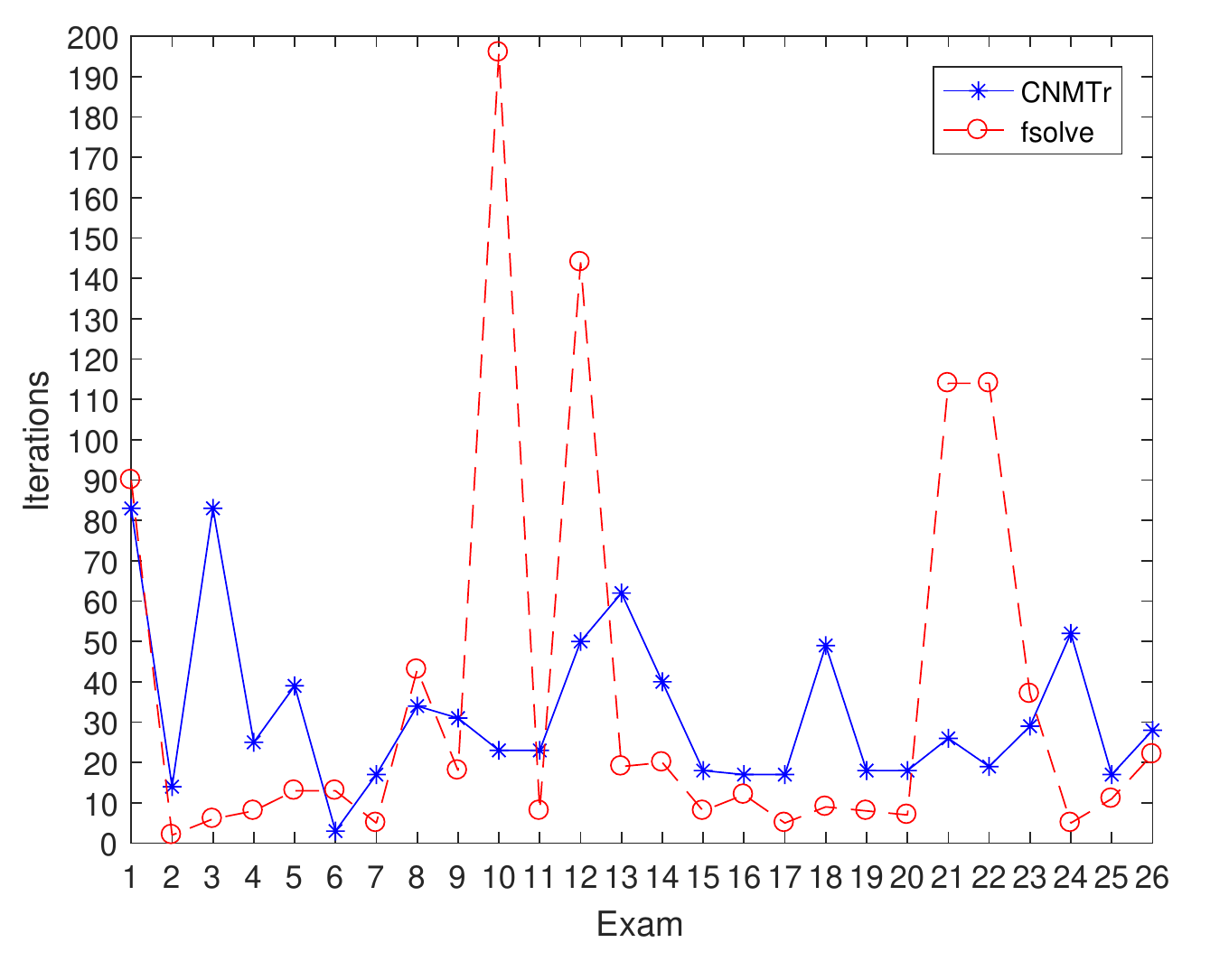}
       \captionof{figure}{The number of iterations.}
       \label{FIGITER}
   \end{minipage}%
   \begin{minipage}{.5\textwidth}
       \centering
       \includegraphics[width=.9\linewidth]{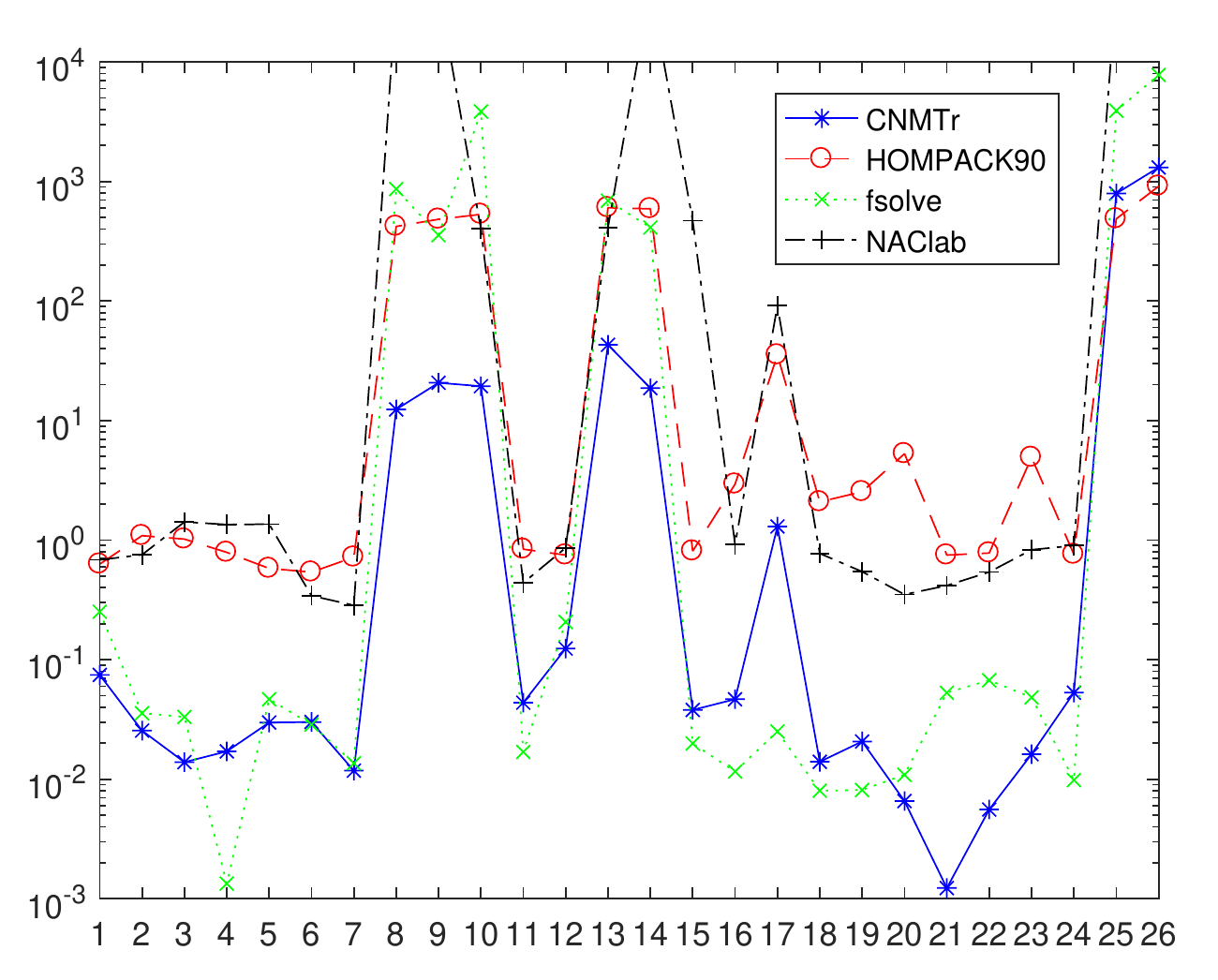}
       \captionof{figure}{The computational time.}
       \label{FIGCONTIME}
   \end{minipage}
\end{figure}

The numerical results are arranged in Table \ref{TABCOM} and Table \ref{TABCOMRE}.
The number of iterations of CNMTr and fsolve is illustrated by Figure \ref{FIGITER}.
The computational time of these four methods (CNMTr, HOMPACK90, fsolve and NAClab)
is illustrated by Figure \ref{FIGCONTIME}. From Table \ref{TABCOM} and Table
\ref{TABCOMRE}, we find that CNMTr performs well for those test problems. However, 
the trust-region method (fsolve) and the classical homotopy methods (HOMPACK90
and NAClab) fail to solve some problems, which especially come from
the real-world problems with the non-isolated singular Jacobian matrices such
as examples $1, \, 2, \, 3, \, 4, \,  6, \, 21, \,  23$. Furthermore, from Figures \ref{FIGITER}
and \ref{FIGCONTIME}, we also find that CNMATr has the same fast convergence property as the 
traditional optimization method (fsolve).
\begin{table}
   \newcommand{\tabincell}[2]{\begin{tabular}{@{}#1@{}}#2\end{tabular}}
  \centering
  \caption{Test problems.}
  \label{TABPROB}
    \begin{tabular}{|c|c|c|}
       \hline
       Problems & dimension & problem descriptions   \\ \hline
       Exam 1 &  $n = 3$  & Robertson problem, an autocatalytic reaction
       \cite{HW1996,Robertson1966} \\ \hline
       Exam 2 & $n = 4$ & E5, the chemical pyrolysis \cite{HW1996} \\ \hline
       Exam 3 & $n = 20$ & The pollution problem \cite{VL1994} \\ \hline
       Exam 4 & $n = 5$ & The stability problem of an aircraft (p. 279, \cite{NW1999}) \\ \hline
       Exam 5 & $n = 1$ & $F(x) = \sin(5x)-x$ (p. 279, \cite{NW1999}) \\ \hline
       Exam 6 & $n = 2$ & \tabincell{c}{$e^{x^2 + y^2} - 3 = 0$, \\
                         $x + y - \sin \big(3(x+y)\big) = 0$ (p. 149, \cite{Deuflhard2004})} \\ \hline
       Exam 7 & $n = 2$ & $x = 0, \; -2y = 0$      \\ \hline
       Exam 8 & $n = 3000$ & Extended Rosenbrock function (p. 362, \cite{DS2009}
       or \cite{MGH1981}) \\ \hline
       Exam 9 & $n = 3000$ & Extended Powell singular function (p. 362, \cite{DS2009}
       or \cite{MGH1981}) \\ \hline
       Exam 10 & $n = 3000$ & Trigonometric function (p. 362, \cite{DS2009}
       or \cite{MGH1981}) \\ \hline
       Exam 11 & $n = 3$ & Helical valley function (p. 362, \cite{DS2009}) \\ \hline
       Exam 12 & $n = 4$ & Wood function (p. 362, \cite{DS2009}) \\ \hline
       Exam 13 & $n = 3000$ & Extended Cragg and Levy function \cite{LUL1994} \\ \hline
       Exam 14 & $n = 3000$ & Singular Broyden problem \cite{LUL1994} \\ \hline
       Exam 15 & $n = 10 $ & The tridiagonal system \cite{LUL1994} \\ \hline
       Exam 16 & $n = 10$ & The discrete boundary-value problem \cite{LUL1994} \\ \hline
       Exam 17 & $n = 100$ & Broyden tridiagonal problem \cite{LUL1994} \\ \hline
       Exam 18 & $n = 5$ & The asymptotic boundary value problem \cite{PHP1975} \\ \hline
       Exam 19 & $n = 3$ & The box problem \cite{MGH1981} \\ \hline
       Exam 20 & $n = 2$ & \tabincell{c}{$f_{1}(x) = x_{1}^{2} + x_{2}^{2} - 2$, \\
                 $f_{2}(x) = e^{x_{1}-1} + x_{2}^{2} - 2$ (p.149, \cite{DS2009})} \\ \hline
       Exam 21 & $n = 2$ & Powell badly scaled function \cite{MGH1981} \\ \hline
       Exam 22 & $n = 2$ & Chemical equilibrium problem 1 \cite{KLH1982} \\ \hline
       Exam 23 & $n = 6$ & Chemical equilibrium problem 2 \cite{KLH1982} \\ \hline
       Exam 24 & $n = 10$ & Brown almost linear function \cite{MGH1981} \\ \hline
       Exam 25 & $n = 3000$ & \tabincell{c}{$a = 2*\text{ones}(n, 1), \; b = \text{ones}((n-1), 1)$, \\
                              $A = \text{diag}(a,1) + \text{diag}(b,1) + \text{diag}(b,-1)$, \\
                              $Ax - \lambda x = 0, \; x^{T}x = 1$}  \\ \hline
       Exam 26 & $n = 3000$ & \tabincell{c}{$a = \text{ones}(n, 1), \; b = \text{ones}((n-1), 1), \; c = 2*b$, \\
                              $A = \text{diag}(a,1) + \text{diag}(b,1) + \text{diag}(c,-1)$, \\
                              $Ax - \lambda x = 0, \; x^{T}x = 1$}  \\ \hline
    \end{tabular}
\end{table}

\renewcommand{\arraystretch}{1.1}
\begin{table}
  \newcommand{\tabincell}[2]{\begin{tabular}{@{}#1@{}}#2\end{tabular}}
  \centering
  \scriptsize
  \caption{Numerical results.}
  \label{TABCOMRE}
    \begin{tabular}{|c|c|c|c|c|c|c|c|c|}
    \hline
    \multirow{2}{*}{Exam}&
    \multicolumn{2}{c|}{CNMTr}&\multicolumn{2}{c|}{HOMPACK90}&\multicolumn{2}{c|}{fsolve}
    &\multicolumn{2}{c|}{NAClab (psolve)}\cr \cline{2-9}
    & CPU (s) & $\|F(x^{it})\|_{\infty}$ & CPU (s) & $\|F(x^{it})\|_{\infty}$ & CPU (s)
    & $\|F(x^{it})\|_{\infty}$ & CPU (s) & $\|F(x^{it})\|_{\infty}$ \cr
    \hline
     1 &7.46E-02& 4.87E-13&6.31E-01& \tabincell{c}{5.01E-04 \\ \red{(failed)}} &2.52E-01&1.64E-07&6.84E-01
     & \tabincell{c}{0 \\ \red{(failed)}} \cr \hline      
     2 &2.55E-02& 9.86E-14 &1.09& \tabincell{c}{3.06 \\ \red{(failed)}}  &3.57E-02
     & \tabincell{c} {1.39E-12 \\ \red{(far sol.)}} &7.55E-01& \tabincell{c}{4.27E-20 \\ \red{(failed)}}
     \cr \hline      
     3 &1.39E-02& 9.61E-14&1.02& \tabincell{c}{3.12 \\ \red{(failed)}} &3.32E-02
     & \tabincell{c}{9.24E-05 \\ \red{(far sol.)}} &1.42& \tabincell{c}{0 \\ \red{(failed)}} \cr\hline      
     4 &1.71E-02& 1.17E-15& 7.94E-01& \tabincell{c}{0.74 \\ \red{(failed)}} &1.34E-03
     & \tabincell{c}{1.93 \\ \red{(failed)}} &1.35& \tabincell{c} {1.40E+01 \\ \red{(failed)}} \cr \hline      			
     5 &2.98E-02& 4.66E-15&5.87E-01&2.60E-12&4.67E-02& \tabincell{c}{5.51E-01 \\ \red{(failed)}} 
     &1.36& \tabincell{c}{1.96 \\ \red{(failed)}} \cr \hline     				 		
    6 &3.01E-02& 1.11E-16&5.49E-01& \tabincell{c} {1.34E-02 \\ \red{(failed)}} &2.87E-02&3.44E-15&3.43E-01
      & \tabincell{c}{4.39 \\ \red{(failed)}} \cr \hline      	 		
    7 &1.18E-02& 3.05E-13&7.52E-01&0&1.36E-02&2.34E-09&2.85E-01&0 \\ \hline
    8 &1.24E+01 & 3.91E-13&4.21E+02&5.12E-13&8.64E+02&3.20E-13&3.65E+04&7.12E-12\\ \hline
    9 &2.07E+01& 4.10E-13&4.83E+02&6.84E-12&3.55E+02&7.50E-13&3.97E+04&5.21E-13 \\ \hline
    10 &1.94E+01& 4.05E-13&5.31E+02&6.31E-15&3.85E+03& \tabincell{c}{6.75 \\ \red{(failed)}} & 4.02E+03 
    & \tabincell{c}{1.20E+04 \\ \red{(failed)}} \cr \hline   
    11 &4.37E-02& 2.58E-13&8.37E-01&2.15E-14&1.69E-02&1.39E-17&4.39E-01
    & \tabincell{c}{9.90E+02 \\ \red{(failed)}} \cr\hline
    12 &1.24E-01& 6.77E-13&7.53E-01&8.94E-13&2.07E-01& \tabincell{c}{5.25E-01 \\ \red{(failed)}} 
    &8.62E-01&6.02E-12 \cr\hline
    13 &4.30E+01& 9.58E-13&6.03E+02&9.68E-13&6.91E+02& \tabincell{c}{4.57E-01 \\ \red{(failed)}} 
    &4.13E+03& \tabincell{c}{4.84E+08 \\ \red{(failed)}} \cr\hline
    14 &1.87E+01& 6.31E-13&5.91E+02&8.41E-13&4.12E+02&1.48E-06&4.10E+04&8.13E-12\\ \hline
    15 &3.80E-02& 1.42E-14&8.06E-01&3.84E-14&1.99E-02&5.72E-13&4.71E+02&5.18E-13\\ \hline
    16 &4.67E-02& 2.44E-14&2.94&6.57E-14&1.16E-02&6.76E-13&9.20E-01&4.15E-13\\ \hline
    17 &1.30& 6.17E-13&3.54E+01&5.71E-13&2.51E-02&8.88E-16&9.13E+01&3.16E-12\\ \hline
    18 &1.40E-02& 3.81E-16&2.09&6.14E-16&8.02E-03&7.19E-12&7.74E-01
    & \tabincell{c}{1 \\ \red{(failed)}} \cr\hline
    19 &2.07E-02& 5.84E-15&2.54&6.58E-12&8.13E-03&4.96E-13&5.45E-01&0 \\ \hline
    20 &6.57E-03& 2.66E-15&5.28&5.14E-13&1.09E-02&2.22E-16&3.51E-01
    & \tabincell{c}{8.53 \\ \red{(failed)}}\cr\hline
    21 &1.23E-03& 8.77E-15&7.53E-01& \tabincell{c}{1.22E-02 \\ \red{(failed)}} &5.29E-02&3.55E-05
    &4.17E-01& \tabincell{c}{1 \\ \red{(failed)}} \cr\hline
    22 &5.58E-03& 0&7.79E-01&0&6.73E-02& \tabincell{c}{2.73 \\ \red{(failed)}} 
    &5.41E-01&0 \cr\hline							
    23 &1.62E-02&4.48E-13&4.92& \tabincell{c}{1.09E+02 \\ \red{(failed)}} &4.84E-02
    & \tabincell{c}{1.05E+02 \\ \red{(failed)}} &8.31E-01& \tabincell{c}{5.47E+14 \\ \red{(failed)}} \cr\hline
    24 &5.30E-02& 1.24E-14&7.63E-01&6.26E-13&9.83E-03&4.23E-12&9.00E-01&2.22E-16\\ \hline
    25 &7.94E+02& 6.09E-13 &4.87E+02&4.13E-13&3.91E+03&3.55E-15&4.07E+04&1.45E-12\\ \hline
    26 &1.31E+03& 2.30E-13&9.12E+02&6.14E-12&7.77E+03&5.55E-16&7.09E+04&6.14E-12\\ \hline
    \end{tabular}
\end{table}

\begin{table}[!htbp]
   \centering
   \caption{Statistical results.}
   \label{TABCOM}
   \begin{tabular}{ccccc}
     \hline
     &CNMTr & HOMPACK90 & fsolve & NAClab\\
     \hline
     number of failed problems   & 0 & 7  & 9  & 13  \\
     number of the minimum time    & 19  & 0  & 7  & 0  \\
     \hline
   \end{tabular}
\end{table}

\section{Conclusions}

In this article, we consider the continuation Newton method with the new 
time-stepping scheme (CNMTr) based on the trust-region updating strategy. We 
also analyze its local and global convergence for the nonsingular Jacobian and 
singular Jacobian problems. Finally, for some classical test problems, we compare 
it with the classical homotopy methods (HOMPACK90 and psolve.m) and the 
traditional optimization method (fsolve.m). Numerical results show that CNMTr 
is more robust and faster than the traditional optimization method. From our 
point of view, the continuation Newton method with the trust-region updating
strategy (Algorithm \ref{CNMTR}) is worth investigating further as a special 
continuation method. We have also extended it to the linear programming problem 
\cite{LY2021}, the unconstrained optimization problem \cite{LXLZ2021} and the 
underdetermined system of nonlinear equations \cite{LX2021}. The promising results 
are reported for those problems therein. 

\vskip 2mm

\noindent \textbf{Acknowledgments} \; This work was supported in part by Grant
61876199 from National Natural Science Foundation of China, and Grant YJCB2011003HI
from the Innovation Research Program of Huawei Technologies Co., Ltd.. 
The authors are grateful to two anonymous referees for their helpful comments
and suggestions.

\end{document}